\documentclass{amsart}
\usepackage{amscd,amsmath,amssymb}
\usepackage[mathscr]{eucal}
\usepackage{graphicx}
\newtheorem{thm}{Theorem}[section]

\newtheorem{cor}[thm]{Corollary}
\newtheorem{lem}[thm]{Lemma}

\theoremstyle{definition}
\newtheorem{defn}[thm]{Definition}
\theoremstyle{remark}
\newtheorem{rem}[thm]{Remark}
\newtheorem{ex}[thm]{Example}
\numberwithin{equation}{section}

\newcommand{\abs}[1]{\left\vert#1\right\vert}

\newcommand{\eps}{\varepsilon}

\begin{document}
\title{The curvature of contact structures on 3-manifolds.}
\author{Vladimir Krouglov.}
\maketitle
\begin{abstract}
We study the sectional curvature of plane distributions on 3-manifolds. We show that if the distribution is a contact structure it is easy to manipulate this curvature. As a corollary we obtain that for every transversally oriented contact structure on a closed 3-dimensional manifold $M$ there is a metric, such that the sectional curvature of the contact distribution is equal to $-1$. We also introduce the notion of Gaussian curvature of the plane distribution. For this notion of curvature we get the similar results.
\end{abstract}

\section{Introduction}
The problem of prescribing the curvatures of the manifold is one of the central problems in Riemmanian geometry. That is, given a smooth function can it be realized as a scalar(Ricci or sectional) curvature of some Riemmanian metric on a manifold. The solution of Yamabe problem is the best known result in prescribing the scalar curvature on a manifold(cf. \cite{LP}). There are several results on prescribing the Ricci curvature of a manifold(cf. for example \cite{L4}). It is natural to ask to what extent it is possible to prescribe the sectional curvature of the plane distribution on a 3-manifold. It turns out that this problem is closely connected with the contactness of the distribution. In fact we have the following: \\ \\
\textbf{Theorem 4.1}
\emph{Let $\xi$ be a transversally orientable contact structure on a closed orientable 3-manifold $M$. For any smooth strictly negative function $f$, there is a metric on $M$, such that $f$ is the sectional curvature of $\xi$.} \\

If we impose more topological restrictions on the distribution we can obtain even stronger result: \\ \\
\textbf{Theorem 5.1}
\emph{Let $\xi$ be a transversally orientable contact structure on a closed orientable 3-manifold $M$, such that there is a contact structure $\eta$ transversal to $\xi$. Then for any smooth function $f$, there is a metric on $M$, such that $f$ is a sectional curvature of $\xi$.} \\

It is a well known problem, whether a foliation on a 3-dimensional manifold admits a simultaneous uniformization of all its leaves. Reeb stability theorem asserts that on a compact orientable 3-manifold the only foliation with the leaves having positive Gaussian curvature is the foliation of $M = S^2 \times S^1$ by spheres. It is known \cite{Can} that if $M$ is atoroidal and aspherical and the foliation is taut, then there is a metric on $M$, such that all leaves have constant negative Gaussian curvature $-1$. In the case of contact structures we ask the similar question. For this we have to  introduce the notion of Gaussian curvature of the plane distribution.

 We define the Gaussian curvature of the plane distribution as the sum $K_G(\xi) = K(\xi) + K_e(\xi)$ of the sectional and the extrinsic curvatures of the distribution. In the case of integrable $\xi$ this equation is nothing else as a Gauss equation.
\begin{defn}
Let $\xi$ be a plane distribution on $M$. We say that $\xi$ admits an uniformization if there is a metric on $M$, such that the Gaussian curvature of $\xi$ is constant.
\end{defn}
It turns out that unlike the case of foliations every transversally orientable contact structure on a closed 3-manifold admits an uniformization. We have the following \\ \\
\textbf{Theorem 6.1}
\emph{Let $\xi$ be a transversally orientable contact structure on a closed orientable 3-manifold $M$. For any smooth strictly negative function $f$, there is a metric on $M$, such that $f$ is the Gaussian curvature of $\xi$.} \\ \\
This paper is organized as follows. In Section 2. we recall basic facts about the geometry of plane distributions. In Section 3. we prove the main technical lemma. Section 4. is devoted to the proof of Theorem 1.1. In Section 5. we prove Theorem 1.2 and some of its corollaries. We prove Theorem 1.4 in Section 6.
\section{Basic Definitions and Notations.}
Throughout this paper $M$ will be a closed orientable 3-manifold. A distribution on $M$ is a two dimensional subbundle of the tangent bundle of $M$. That is, at each point $p$ in $M$ there is a plane $\xi_p$ in the tangent space $T_p M$. A distribution is called integrable, if there is a foliation on $M$ which is tangent to it. The following Frobenius theorem gives necessary and sufficient conditions for $\xi$ to be integrable.
\begin{thm}
Let $\xi$ be a distribution on $M$. Then $\xi$ is integrable if and only if for any two sections $S$ and $T$ of $\xi$ its Lie bracket belongs to $\xi$.
\end{thm}
\begin{defn}
A distribution $\xi$ is called a contact structure if for any linearly independent sections $S$ and $T$ of $\xi$ and for any $p \in M$ the Lie bracket $[S, T]$ at $p$ does not belong to $\xi_p$.
\end{defn}
A distribution $\xi$ is called transversally oriented if there is a globally defined $1$-form $\alpha$, such that $\xi = Ker(\alpha)$ . This is equivalent to say that there exists a globally defined vector field $n$ which is transverse to $\xi$. It is an easy consequence of Frobenius Theorem that $\xi$ is a contact structure if and only if
$$
\alpha \wedge d\alpha \ne 0
$$
Fix some orientation on $M$. A contact structure is said to be positive(resp. negative) if the orientation induced by $\alpha \wedge d\alpha$ coincides(resp. is opposite to) the orientation on $M$.

A contact structure $\xi$ is called overtwisted, if there is an embedded disk, such that $TD|_{\partial D} = \xi|_{\partial D} $. If $\xi$ is not overtwisted, it is called tight.  

The Euler class $e(\xi) \in H^2(M,\mathbb{Z})$ of the plane distribution is the Euler class of the fibration $\xi \to M$.
It is known that if $\xi$ is a 2-dimensional plane distribution on $M$ with vanishing Euler class then $\xi$ is trivial.
Recall, that a framing of $M$ is the presentation of the tangent bundle of $M$ as a product $TM \simeq M \times \mathbb{R}^3$.
A framing on $M$ consists of three linearly independent vector fields. It is known that every closed orientable 3-manifold admits a framing.

A bi-contact structure on $M$ is a pair $(\xi, \eta)$ of transverse contact structures, such that they define opposite orientation on $M$.

Assume that $M$ is a Riemannian manifold with the metric $\langle \cdot, \cdot \rangle$ and the Levi-Chevita connection $\nabla$. Let $n$ be a local unit vector field orthogonal to $\xi$. We are now going to define the second fundamental form of $\xi$. The definition is due to Reinhart \cite{Re}.
\begin{defn}
The second fundamental form of $\xi$ is a symmetric bilinear form, which is defined in the following way:
$$
B(S, T) = \frac{1}{2}\langle \nabla_S T + \nabla_T S, n \rangle
$$
for all sections $S$ and $T$ of $\xi$.
\end{defn}
\begin{rem} If $\xi$ is integrable, then $B$ restricted to the leaf of $\xi$ coincides with the second fundamental form of the leaf.
\end{rem}
Let $S$ and $T$ be two linearly independent sections of $\xi$.
\begin{defn}We call the function
$$
K_e(\xi) = \frac{B(S, S)B(T, T) - B(S, T)^2}{\langle S,S\rangle \langle T,T\rangle - \langle S, T\rangle^2}
$$ an extrinsic curvature of $\xi$.
\end{defn}
It is easy to verify that $K_e(\xi)$ depends only on $\xi$, not on the actual choice of $S$, $T$ and $n$.
\begin{defn}Consider the function $K(\xi)$ which assigns to a point $p \in M$ the sectional curvature of the plane $\xi_p$. We call this function the sectional curvature of $\xi$.
\end{defn}
\begin{defn}
We call the sum $K_G(\xi) = K(\xi) + K_e(\xi)$ the Gaussian curvature of $\xi$.
\end{defn}
Let $S$, $T$ and $U$ be the local sections of $TM$. Recall the Koszul formula for the Levi-Chevita connection of $\langle \cdot, \cdot \rangle$:
$$
2 \langle\nabla_S T, U\rangle = S\langle T, U\rangle + T\langle U,
S \rangle - U\langle S, T\rangle + \langle[S, T], U\rangle-$$
$$ -\langle [S,U], T\rangle - \langle[T,U], S\rangle
$$
\section{The Deformation of Metric.}
In this section we will give the proof of the main technical results we will need throughout the paper.

Let $\xi$ be a transversally orientable plane distribution on a 3-dimensional Riemmanian manifold $(M, \langle\cdot, \cdot \rangle)$. Fix a unit normal vector field $n$. Suppose $a$ is a strictly positive smooth function on $M$. A stretching of $\langle \cdot, \cdot \rangle$ along $n$ by the function $a$ is the following Riemmanian metric on $M$:
$$
\langle \cdot, \cdot \rangle_a=a \langle \cdot, \cdot \rangle|_n \oplus \langle\cdot, \cdot \rangle|_\xi
$$
Our aim is to calculate the sectional curvature of $\xi$ in the stretched metric
in terms of the initial metric.

Consider an open subset $U \subset M$, such that $\xi|_U$ is a trivial fibration. Let $X$ and $Y$ be a pair of orthonormal sections of $\xi|_U$. The triple $(X, Y, n)$ is an orthonormal framing on $U$ with respect to $\langle \cdot, \cdot \rangle$.

In the stretched metric this frame is orthogonal, vector fields $X$ and $Y$ are unit and the length of $n$ is equal to $a$.
Denote by $\nabla$ the Levi-Chevita connection of $\langle \cdot, \cdot \rangle_a$.
\begin{lem} The sectional curvature of $\xi$ with respect to $\langle\cdot, \cdot \rangle_a$ can be calculated by the following formula:
$$
K(\xi) = - \frac{3}{4}a \langle[X, Y], n\rangle^2 + (- X
\langle[X, Y], Y \rangle - Y \langle[X, Y], X \rangle - \langle
[X,Y], X\rangle^2 - \langle [X,Y], Y\rangle^2 + $$
$$ + \frac{1}{2}
\langle [X,Y], n\rangle (-\langle [n, Y], X \rangle + \langle [n,
X], Y \rangle)) +
$$
$$ + \frac{1}{a} (\frac{1}{4}(\langle [X, n], Y
\rangle + \langle [Y, n], X \rangle)^2 - \langle [Y, n], Y\rangle
\langle[X, n], X\rangle)
$$
\end{lem}
\emph{Proof:} Since $X$ and $Y$ are unit, the sectional curvature of $\xi$ is calculated by the formula
$$
K(\xi) =\langle R(X, Y)Y, X\rangle_a = \langle \nabla_X \nabla_Y
Y, X \rangle_a - \langle \nabla_Y \nabla_X Y, X\rangle_a - \langle
\nabla_{[X, Y]} Y, X\rangle_a
$$
The first summand can be rewritten:
$$
\langle \nabla_X \nabla_Y Y, X \rangle_a = X \langle \nabla_Y Y,
X\rangle_a - \langle \nabla_Y Y, \nabla_X X \rangle_a =
$$
Apply the Koszul formula to $X \langle \nabla_Y Y, X\rangle_a$. We get
$$
  = X\langle
[Y, X], Y\rangle_a - \langle \langle \nabla_Y Y, n \rangle_a \frac{n}{a} + \langle \nabla_Y Y,
 Y\rangle_a Y
  + \langle \nabla_Y Y, X\rangle_a X, \nabla_X X \rangle_a =
$$
Since $X$ and $Y$ are of unit length this reduces to:
$$
= - X \langle[Y, X], Y \rangle_a - \frac{1}{a}\langle \nabla_Y Y,
n \rangle_a\langle \nabla_X X, n\rangle_a =
$$
Apply the Koszul formula to $\langle \nabla_Y Y,
n \rangle_a\langle \nabla_X X, n\rangle_a$. This gives us:
$$
= - X \langle[Y, X], Y
\rangle_a - \frac{1}{a}\langle [Y, n], Y\rangle_a \langle[X, n],
X\rangle_a
= - X \langle[Y, X], Y
\rangle - \frac{1}{a}\langle [Y, n], Y\rangle \langle[X, n],
X\rangle
$$
The second summand is equal to:
$$
 - \langle \nabla_Y \nabla_X Y, X \rangle_a = - Y \langle \nabla_X Y,
X\rangle_a + \langle \nabla_X Y, \nabla_Y X \rangle_a =  Y\langle
Y, \nabla_X X\rangle_a +
$$
$$
 + \langle \langle \nabla_X Y, n \rangle_a \frac{n}{a} + \langle \nabla_X Y,
 Y\rangle_a Y
  + \langle \nabla_X Y, X\rangle_a X, \nabla_Y X \rangle_a =
$$
$$
= - Y \langle[X, Y], X \rangle_a + \frac{1}{a}\langle \nabla_X Y, n
\rangle_a \langle \nabla_Y X, n\rangle_a =
$$
Write the equations for the derivatives $\langle \nabla_X Y
, n \rangle_a$ and $\langle \nabla_Y X , n \rangle_a$:
$$
2 \langle \nabla_X Y , n \rangle_a = \langle [X, Y], n \rangle_a -
\langle [X, n], Y \rangle_a - \langle [Y, n], X \rangle_a =
$$
$$ a \langle [X, Y], n \rangle - \langle [X, n], Y \rangle -
\langle [Y, n], X \rangle
$$
$$
2 \langle \nabla_Y X , n \rangle_a = \langle [Y, X], n \rangle_a -
\langle [Y, n], X \rangle_a - \langle [X, n], Y \rangle_a =
$$
$$
 a \langle [Y, X], n \rangle - \langle [Y, n], X \rangle -
\langle [X, n], Y \rangle
$$
Inserting the above equations into the second summand we have:
$$
=  - Y \langle[X, Y], X \rangle_a + \frac{1}{4a}( - a \langle [X, Y],
n \rangle + \langle [X, n], Y \rangle + \langle [Y, n], X \rangle)
( - a \langle [Y, X], n \rangle +
$$
$$
+\langle [Y, n], X \rangle +
\langle [X, n], Y \rangle)
$$
The last summand is:
$$
- \langle \nabla_{[X, Y]} Y, X\rangle_a = - \langle \nabla_{\langle
[X,Y], n\rangle n + \langle [X,Y], X\rangle X + \langle [X,Y],
Y\rangle Y } Y, X \rangle_a =
$$
$$
 = - \langle
[X,Y], n\rangle \langle \nabla_n Y, X\rangle_a  - \langle [X,Y],
X\rangle \langle \nabla_X Y, X \rangle_a - \langle [X,Y], Y\rangle
\langle \nabla_Y Y, X\rangle_a
$$
The derivative $\langle \nabla_n Y,
X\rangle_a$ is equal to:
$$
\langle \nabla_n Y, X\rangle_a = - \frac{1}{2}( - \langle [n, Y],
X \rangle_a + \langle [n, X], Y \rangle_a + \langle [Y, X], n
\rangle_a) =
$$
$$
= - \frac{1}{2}( - \langle [n, Y], X \rangle + \langle [n, X], Y
\rangle + a\langle [Y, X], n \rangle)
$$
Which gives us:
$$
- \langle \nabla_{[X, Y]} Y, X\rangle_a = - \langle [X,Y], n\rangle
\langle \nabla_n Y, X\rangle_a  - \langle [X,Y], X\rangle \langle
\nabla_X Y, X \rangle_a -
$$
$$
- \langle [X,Y], Y\rangle \langle \nabla_Y
Y, X\rangle_a
= \frac{1}{2} \langle [X,Y], n\rangle( - \langle [n, Y], X
\rangle + \langle [n, X], Y \rangle + a\langle [Y, X], n \rangle) -
$$
$$
- \langle [X,Y], X\rangle^2 - \langle [X,Y], Y\rangle^2
$$
Summing this up ,the sectional curvature of $\xi$ is equal to:
$$
K(\xi) = - X \langle[Y, X], Y \rangle - \frac{1}{a}\langle [Y,
n], Y\rangle \langle[X, n], X\rangle
- (Y \langle[X, Y], X \rangle - $$$$ - \frac{1}{4a}(-a \langle [X, Y], n
\rangle + \langle [X, n], Y \rangle + \langle [Y, n], X \rangle)
(-a \langle [Y, X], n \rangle + \langle [Y, n], X \rangle +
\langle [X, n], Y \rangle)) -
$$
$$
 (- \frac{1}{2} \langle [X,Y], n\rangle(- \langle [n, Y], X \rangle
+ \langle [n, X], Y \rangle + a\langle [Y, X], n \rangle) +
\langle [X,Y], X\rangle^2 + \langle [X,Y], Y\rangle^2)
$$
It is straightforward to verify that this gives us the desired expression.

\begin{lem}
The extrinsic curvature $K_e(\xi)$ with respect to $\langle\cdot, \cdot \rangle_a$ can be calculated by the following formula:
$$
K_e(\xi) =\frac{1}{a}( \langle[X, n], X\rangle \langle[Y, n],
Y\rangle - \frac{1}{4}(\langle[X, n], Y\rangle + \langle[Y, n],
X\rangle)^2)
$$

\end{lem}
\emph{Proof:} Since $X$ and $Y$ are unit the extrinsic curvature
$$K_e(\xi) = B(X, X)B(Y, Y) - B(X, Y)^2$$. By the definition of $B$ the extrinsic curvature is equal to:
$$
K_e(\xi) = \langle \nabla_X X , \frac{n}{\sqrt{a}} \rangle_a \langle\nabla_Y Y, \frac{n}{\sqrt{a}} \rangle_a - \frac{1}{4} \langle \nabla_X Y + \nabla_Y X, \frac{n}{\sqrt{a}} \rangle_a^2
$$
Apply the Koszul formula to covariant derivatives:
$$
K_e(\xi) =\frac{1}{a}( \langle[X, n], X\rangle_a \langle[Y, n],
Y\rangle_a - \frac{1}{4}(\frac{1}{2} \langle[X, Y], n\rangle_a - \frac{1}{2}\langle [X, n], Y \rangle_a -
\frac{1}{2} \langle [Y, n], X \rangle_a -
$$
$$
 - \frac{1}{2}\langle [Y, X], n\rangle_a
- \frac{1}{2}\langle[Y, n], X\rangle_a - \frac{1}{2} \langle [X, n], Y \rangle_a)^2 =
$$
$$
= \frac{1}{a}( \langle[X, n], X\rangle \langle[Y, n],
Y\rangle - \frac{1}{4}(\langle [X, n], Y \rangle + \langle [Y, n], X \rangle)^2
$$

Summing the extrinsic curvature of $\xi$ with the sectional curvature gives us the Gaussian curvature of the plane distribution $\xi$.
\begin{lem}The Gaussian curvature $K_G(\xi)$ can be calculated by the formula:
$$
K_G(\xi) =  K(\xi) + K_e(\xi) = - \frac{3}{4}a \langle[X, Y], n\rangle^2 +
(X \langle[X, Y], Y \rangle - Y \langle[X, Y], X \rangle - \langle
[X,Y], X\rangle^2 -
$$$$
 -\langle [X,Y], Y\rangle^2)
+
\frac{1}{2} \langle [X,Y], n\rangle (-\langle [n, Y], X \rangle +
\langle [n, X], Y \rangle))
$$
\end{lem}
\begin{rem}
If $\xi$ is integrable then $\langle [X, Y], n\rangle = 0$ and
$$
K_G(\xi) = (X \langle[X, Y], Y \rangle - Y \langle[X, Y], X \rangle - \langle
[X,Y], X\rangle^2 -\langle [X,Y], Y\rangle^2
$$ is nothing else as the expression of the Gaussian curvature of the leaves of $\xi$ written in the local frame tangent to the leaves.
\end{rem}
\begin{lem}
Let $(X, Y, n)$ be a framing on $M$. Assume that distribution spanned by $n$ and $Y$ is a contact structure. Then there is a metric on $M$, such that extrinsic curvature of the distribution spanned by $X$ and $Y$ is strictly less then zero.
\end{lem}
\emph{Proof:} Fix a metric $\langle \cdot, \cdot \rangle$ such that the framing is orthonormal. Let $\xi$ be a distribution spanned by vector fields $X$ and $Y$. Stretch the metric along $X$ by a constant factor $\lambda^2$ and along $Y$ by a constant factor $\frac{1}{\lambda^2}$. Calculate the extrinsic curvature of $\xi$ in the stretched metric $\langle \cdot, \cdot \rangle_\lambda$
$$
K_e(\eta) = \langle[n, X], X\rangle_\lambda \langle[n, Y],
Y\rangle_\lambda - \frac{1}{4}(\langle[n, X], Y\rangle_\lambda +
\langle[n, Y], X\rangle_\lambda)^2=
$$
$$
= \lambda^2 \langle[n, X], X\rangle \frac{1}{\lambda^2}\langle[n,
Y], Y\rangle - \frac{1}{4}(\frac{1}{\lambda^2}\langle[n, X],
Y\rangle + \lambda^2 \langle[n, Y], X\rangle)^2 =
$$
$$
=\langle[n, X], X\rangle\langle[n, Y], Y\rangle -
\frac{1}{4}(\frac{1}{\lambda^2}\langle[n, X], Y\rangle + \lambda^2
\langle[n, Y], X\rangle)^2= $$$$ =\langle[n, X],
X\rangle\langle[n, Y], Y\rangle -
\frac{1}{4\lambda^4}\langle[n,X], Y\rangle^2 -
\frac{1}{2}\langle[n, X], Y\rangle\langle[n, Y], X\rangle-
\frac{\lambda^4}{4}\langle[n, Y], X\rangle^2
$$

Since $M$ is compact there is a positive constant $C$ such that
$$
\abs{\langle[n, X], X\rangle \langle[n, Y], Y\rangle -
\frac{1}{2}\langle[n, X], Y\rangle\langle[n, Y], X\rangle} < C
$$
We assumed that distribution spanned by vector fields $n$ and $Y$ is a contact structure. The form $\alpha(\ast) = \langle \ast, X \rangle$ is a contact form of this distribution, so $\langle[n, Y], X\rangle = \alpha([n, Y])
\ne 0$. Since $M$ is compact there is an $\eps$ such that
$$
\abs{\langle[n, Y], X\rangle} > \eps
$$
This means that
$$
K_e(\eta) < C - \frac {\lambda^4\eps^2}{4}
$$
This expression is strictly negative for some sufficiently large $\lambda$.
\section{Prescribing the Sectional Curvature of $\xi$.}
\begin{thm}
Let $\xi$ be a transversally orientable contact structure on a closed orientable 3-manifold $M$. For any smooth strictly negative function $f$, there is a metric on $M$, such that $f$ is the sectional curvature of $\xi$.
\end{thm}
\emph{Proof:}
Since $\xi$ is transversally orientable, there is a globally defined vector field $n$ which is transverse to $\xi$. Fix some Riemmanian metric $\langle \cdot, \cdot \rangle$ on $M$, such that $n$ is a unit normal vector field.
Since $M$ is compact, there is a following finite covering of $M$ by the open sets:

1. $M = \bigcup_\alpha U_\alpha$.

2. Each $U_\alpha$ has a compact closure $V_\alpha$.

3. For each $\alpha$ there is an open set $U_\alpha'$, such that $V_\alpha \subset U_\alpha'$ and $\xi_{U_\alpha'}$ is a trivial fibration.

In each $U_\alpha'$ choose an orthonormal framing $(X_\alpha, Y_\alpha, n|_{U_\alpha'})$. Consider the stretching $\langle \cdot, \cdot \rangle_a$ of $\langle \cdot, \cdot \rangle$ along $n$ by a positive function $a$.

According to Lemma $1.1$ the sectional curvature $K(\xi)$ on $U_\alpha'$ can be rewritten in the following way:
$$
K(\xi) = -\frac{3}{4}a\langle [X_\alpha, Y_\alpha], n\rangle^2 + P_\alpha + \frac{1}{a}Q_\alpha
$$
where  $P_\alpha$ and $Q_\alpha$ are the functions on $U_\alpha'$ independent of $a$.

Since $\xi$ is a contact structure and $U_\alpha$ has a compact closure the expression $\langle [X_\alpha, Y_\alpha], n\rangle^2$ is bounded below by some positive constant $\eps$ and the functions $P_\alpha$ and $Q_\alpha$ are bounded from above. Therefore there is a sufficiently large $D_\alpha$, such that the equation
$$
 -\frac{3}{4}a \langle [X_\alpha, Y_\alpha], n\rangle^2 + P_\alpha + \frac{1}{a}Q_\alpha = fD_\alpha
$$
has a strictly positive solution $a_{ \alpha}(D_\alpha)$. Notice, that for any $D>D_\alpha$ this equation still has a positive solution $a_{\alpha}(D)$. Let $D_0 = max_\alpha\{D_\alpha\}$. Solve the equation above for $D_0$ in each chart $U_\alpha$. Let $a_\alpha = a_{\alpha}(D_0)$.

We claim that $a_\alpha$ constructed this way does not depend on the choice of the orthonormal framing $(X_\alpha, Y_\alpha, n|_{U_\alpha})$.
Let $(X^\prime_{\alpha}, Y^\prime_{\alpha}, n|_{U_\alpha})$ be any other orthonormal framing on $\xi|_{U_\alpha}$. This defines a map
$$
\phi_\alpha : U_\alpha \to O(2)
$$
which maps a point $p \in U_\alpha$ to the transition matrix $A_\alpha(p)$ between two framings $(X^\prime_{\alpha}, Y^\prime_{\alpha})$ and $(X_\alpha, Y_\alpha)$ on $\xi$.
It is easy to verify that
$$
\langle [X^\prime_\alpha, Y^\prime_\alpha], n\rangle^2 = \langle [A_\alpha (X_\alpha), A_\alpha (Y_\alpha)], n\rangle^2 = det A_\alpha^2\langle [X_\alpha, Y_\alpha], n\rangle^2 = \langle [X_\alpha, Y_\alpha], n\rangle^2
$$ hence
does not depend on the choice of orthonormal framing. The expression $\frac{1}{a}Q_\alpha = -  K_e(\xi)$ also does not depend on the choice of the trivialization. Finally the sectional curvature $K(\xi)$ is independent of the framing. It is obvious that the right hand side of
$$
P_\alpha = K(\xi) - \frac{1}{a}Q_\alpha + \frac{3}{4}a\langle [X_\alpha, Y_\alpha], n\rangle^2
$$
does not depend on the choice of framing, so does $P_\alpha$.

Therefore, the functions $a_\alpha$ agree on the overlaps and define a global function $a$ on $M$. The sectional curvature of $\xi$ in the metric $\langle \cdot, \cdot \rangle_a$ is $fD_0$. Consider the metric $\langle \cdot, \cdot \rangle_0 = \frac{1}{\sqrt{D_0}}\langle \cdot, \cdot \rangle_a$. It is easy to calculate, that the sectional curvature of $\xi$ in this metric is equal to $f$.
\begin{cor}
For any transversaly orientable contact structure on a closed orientable 3-manifold, there is a metric on $M$,such that the sectional curvature of $\xi$ in this metric is equal to $-1$.
\end{cor}
\section{Plane Distributions Transverse to a Contact Structure.}
\begin{thm}
Let $\xi$ be a transversally orientable contact structure on $M$, such that there is a contact structure $\eta$ transverse to $\xi$. Then for any smooth function $f$, there is a metric on $M$, such that $f$ is a sectional curvature of $\xi$.
\end{thm}
\emph{Proof:} Since $\xi$ is transverse to $\eta$, both $\xi$ and $\eta$ are trivial fibrations. Consider the following vector fields on $M$. Let $X$ be the section of $\xi \bigcap \eta$, $Y$ -- the section of $\xi$ transverse to $\eta$ and $n$ -- the section of $\eta$ transverse to $\xi$. There is a unique metric on $M$, such that the triple $(X, Y, n)$ defines an orthonormal framing. According to Lemma $3.5$ there is a metric $\langle \cdot, \cdot \rangle$ in which the extrinsic curvature of $\xi$ is a strictly negative function.

Consider the stretching of $\langle \cdot, \cdot \rangle$ along $n$ by a positive function $a$. According to Lemma $3.1$, we have to find $a$ to satisfy the equation
$$
-\frac{3}{4}a \langle [X, Y], n \rangle^2 + P - \frac{1}{4a} K_e(\xi) = f
$$
where $P$ is a function on $M$ which is independent of $a$.

But since $- K_e(\xi) > 0$ this equation always has a strictly positive solution $a$. This completes the proof of the theorem.
\begin{rem}
In the proof of Theorem 5.1 it is crucial that $\xi$ itself is a contact structure. At the points where $\langle [X, Y], n \rangle = 0$ the equation may not have any positive solutions.
\end{rem}
The following theorem has been announced in \cite{DAN}:
\begin{thm} \cite{DAN}
Let $\xi$ an $\eta$ be a positive and negative overtwisted contact structures contained in a same homotopy class of plane fields with Euler class zero. Then, we can isotop $\xi$ and $\eta$ so that $(\xi, \eta)$ is a bi-contact structure.
\end{thm}
Combining this result together with Theorem $5.1$ gives us the next
\begin{cor}
For every overtwisted contact structure $\xi$ with Euler class zero, there is a metric on $M$ such that $\xi$ has a constant positive(or zero) sectional curvature.
\end{cor}
\emph{Proof:} From Theorem $5.3$ we can isotop $\xi$ so that it is contained in a bi-contact structure. Apply Theorem $5.1$ and use the pull-back metric $\phi_t^* \langle \cdot, \cdot \rangle_a$ to get the desired metric.

\begin{ex}Propeller construction, \cite{Mit}.
\end{ex}
Consider the following pair of contact structures on $\mathbb{T}^3$
$$
\xi = Ker(\alpha = \cos z dx - \sin z dy + dz)
$$
$$
\eta = Ker(\beta = \cos z dx + \sin z dy)
$$
It is easy to verify, that $\xi$ is transverse to $\eta$ and we get a bi-contact structure. From Theorem $5.1$, there is a metric on $\mathbb{T}^3$ such that $\xi$ has a positive sectional curvature. This is an example of tight contact structure of positive sectional curvature.

\section{Uniformization of Contact Structures on 3-manifolds.}
The same technique as in Theorem $4.1$ can be applied to the Gaussian curvature of contact structures on three manifolds.
\begin{thm}
Let $\xi$ be a transversally orientable contact structure on a closed orientable 3-manifold $M$. For any smooth strictly negative function $f$, there is a metric on $M$, such that $f$ is the Gaussian curvature of $\xi$.
\end{thm}
\emph{Proof:} Same as Theorem 4.1. The only difference is that in the present case we obtain linear equation in $a$ instead of quadratic.
\begin{cor}(Uniformization of Contact Structures.)
For every transversally orientable contact structure $\xi$ on $M$, there is a metric such that $K_G(\xi) = -1$.
\end{cor}
\begin{ex}
Contact structure with $K_G(\xi) = 1$.
\end{ex}
Consider the unit sphere $S^3 \subset \mathbb{C}^2$ with a bi-invariant metric. The standard contact structure on $S^3$ is defined as the kernel of 1-form
$$
\alpha = \sum^2_{i = 1} (x_i dy_i - y_i dx_i)
$$
restricted from $\mathbb{C}^2$ to $S^3$.

This contact structure is orthogonal to a left-invariant vector  field and therefore is left-invariant. Let $(X, Y)$ be a pair of orthonormal left-invariant  sections of $\xi$. Since the metric is bi-invariant
$$
\nabla _S T = \frac{1}{2}[S, T]
$$
for any left-invariant vector fields on $S^3$. Therefore the second fundamental form of $\xi$ vanishes and
$
K_G(\xi) = K(\xi) = 1
$

\end{document}